\documentclass{amsart}
\usepackage{graphicx}
\graphicspath{{ocfigures/}}
\usepackage{amsthm,amsfonts,amssymb,amsmath,picins
}
\usepackage[centerlast]{subfigure}
\usepackage[rotate,arrow,matrix]{xy}

\DeclareMathOperator{\id}{id}

\DeclareMathOperator{\Hom}{Hom}
\DeclareMathOperator{\ori}{or}

\newtheorem{thm}{Theorem}[section]
\newtheorem{conj}[thm]{Conjecture}
\newtheorem{prop}[thm]{Proposition}

\theoremstyle{definition}
\newtheorem{df}[thm]{Definition}

\theoremstyle{remark}

\newtheorem*{rem}{Remark}
\newtheorem{ack}{Acknowledgments}

\def\v#1{\overrightarrow{#1}}
\def\F{\mathcal F}
\def\Z{\mathbb Z}
\def\Q{\mathbb Q}

\def\C{\mathbb C}
\newcommand{\g}{\mathfrak{g}}
\def\M{\mathcal M}
\def\MM{\overline{\M}_{g,b}^{n,\v{m}}}

\def\S{\mathfrak S}
\def\Mg{\M_{g,b}^{n,\v{m}}}
\def\Mgm#1{\M_{g,b}^{n,(#1)}}

\def\MUm{\underline{\M}^{n,m}_b}

\newcommand{\MUmnv}[3]{{\underline{\M}_{#3}^{#1,#2}}}
\def\MUg{\underline{\M}_{g,b}^{n,\v{m}}}
\def\MUgm{\underline{\M}_{g,b}^{n,m}}

\newcommand{\abs}[1]{{\lvert#1\rvert}}
\newcommand{\Cg}{C^{\operatorname{geom}}}
\newcommand{\delco}{\Delta_{\operatorname{co}}}
\newcommand{\End}{{\operatorname{End}}}
\newcommand{\im}{{\operatorname{Im}}}
\newcommand{\del}{\partial}

\begin{document}

\title[Open-Closed Moduli Spaces and Related Algebraic
Structures]{Open-Closed Moduli Spaces and Related Algebraic
  Structures}
\author[E. Harrelson]{Eric Harrelson}
\address {Department of Mathematics\\Stony Brook University\\
Stony Brook, NY 11794}
\email{harrelson@math.sunysb.edu}
\author[A. A. Voronov]{Alexander A. Voronov}
\address {School of Mathematics\\University of Minnesota\\
Minneapolis, MN 55455, USA and Institut des Hautes \'Etudes
Scientifiques\\Le Bois-Marie\\35, route de Chartres\\91440
Bures-sur-Yvette\\France}
\email{voronov@umn.edu}
\author[J. J. Zuniga]{J. Javier Zuniga}
\thanks{Partially supported by a University of Minnesota Doctoral
  Dissertation Fellowship}
\address{Department of Mathematics\\Purdue University\\
West Lafayette, IN 47907}
\email{jzuniga@math.purdue.edu}

\date{June 5, 2008}

\begin{abstract}
  We set up a Batalin-Vilkovisky Quantum Master Equation (QME) for
  open-closed string theory and show that the corresponding moduli
  spaces give rise to a solution, a generating function for their
  fundamental chains. The equation encodes the topological structure
  of the compactification of the moduli space of bordered Riemann
  surfaces. The moduli spaces of bordered $J$-holomorphic curves are
  expected to satisfy the same equation, and from this viewpoint, our
  paper treats the case of the target space equal to a point. We also
  introduce the notion of a symmetric Open-Closed Topological
  Conformal Field Theory (OC TCFT) and study the $L_\infty$ and
  $A_\infty$ algebraic structures associated to it.
\end{abstract}

\maketitle


\section{Introduction}

The string-field theoretic formulation of open-closed string theory
was developed by Zwiebach in \cite{zwiebach:98}. In his seminal work a
suitable BV algebra was introduced and a solution to the QME was
obtained from so-called string vertices. String vertices are cycles in
the infinite dimensional moduli spaces of nonsingular bordered Riemann
surfaces with a choice of parameterization of the boundary. These
cycles must have a certain topological type, basically, the one
governed by the QME, and their construction is far from being
complete. An alternative approach to this theory involves the use of
certain real compactifications of moduli spaces of Riemann surfaces
with boundary. This was achieved by Costello in \cite{costello:gp} for
the closed case using the compactifications introduced in \cite{ksv}.
We obtain an analogous result for the open-closed case in this paper
using the open-closed moduli spaces of \cite{liu}. It is interesting
to note the existence of an interaction between closed and open
strings that does not seem to have appeared before in the literature.
Such interaction arises from the intrinsic nature of degenerations of
surfaces with boundary and is expressed in our work by the component
$\delco$ of the BV operator $\Delta$. In the work \cite{zwiebach:98}
of Zwiebach, this interaction was present implicitly, via the
antibracket with the unstable moduli space of disks with one interior
puncture.

As pointed out by Sullivan in \cite{sullivan:sigma}, his QME satisfied
by the Gromov-Witten potential in the purely closed case is expected
to generalize to an open-closed version of his sigma model. Fukaya has
a solution for the punctured disk case in \cite{fukaya:06}. The goal
of our work is to set out foundations for the study of the full-blown,
arbitrary genus, arbitrary number of closed and open strings
Gromov-Witten theory.

Another accomplishment of this paper is a new treatment of the
algebraic counterpart (the state space) of OC TCFT
(Section~\ref{algebraic}), a new notion of a symmetric OC TCFT, and a
description of algebraic structures associated with an OC TCFT in
Section~\ref{sh}.

\begin{ack}
  We are very grateful to Kevin Costello, Anton Kapustin, Andrey
  Lazarev, Melissa Liu, Andrei Losev, Albert Schwarz, and Scott Wilson
  for helpful discussions. The third author also thanks IH\'ES for
  hospitality during the final stage of writing the paper.
\end{ack}

\subsection{Convention on Chains}
\label{chains}

Throughout the paper we will use a notion of geometric chains $(P,f)$
based on continuous maps $f: P \to X$ from oriented smooth orbifolds
$P$ with corners to a given topological space $X$. The reason for the
use of this notion is technical: we need the unit circle $S^1$ to have
a canonical fundamental cycle and the moduli spaces which we consider
to have canonical fundamental chains, irrespective of the choice of a
triangulation. On the other hand, we believe it is generally better to
work at the more basic level of chains rather than that of homology.
Similar, though different theories have been used by Gromov
\cite{gromov:filling}, Sen and Zwiebach
\cite{zwiebach:93,sen-zwiebach,zwiebach:98}, Fukaya, Oh, Ohta, and Ono
\cite{fukaya:96,fooo}, Jakob \cite{jakob:alt}, Chataur \cite{chataur},
and Sullivan \cite{sullivan:sigma}. Our notion of chains leads to a
version of oriented bordism theory via passing to homology.  If we
impose extra equivalence relations, such as some kind of a suspension
isomorphism, as in \cite{jakob:alt}, or work with piecewise smooth
geometric chains and treat them as currents, in the spirit of
\cite{fooo}, we may obtain a complex whose homology is isomorphic to
the ordinary real homology of $X$.

Given a topological space $X$, a (\emph{geometric}) \emph{chain} is
a formal linear combination over $\Q$ of continuous maps
\[
f: P \to X,
\]
where $P$ is a compact connected oriented (smooth) orbifold with
corners, modulo the equivalence relation induced by isomorphisms
between the source orbifolds $P$. Here, an \emph{orientation} on an
orbifold with corners is a trivialization of the determinant of its
tangent bundle. Geometric chains form a graded $\Q$-vector space
$\Cg_\bullet (X)$, graded by the dimension of $P$.  The boundary of a
chain is given by $(\del P, f|_{\del P})$, where $\del P$ is the sum
of codimension one faces of $P$ with the induced orientation (Locally,
in positively oriented coordinates near $\del P$, the manifold $P$ is
given by the equation ``the last coordinate is nonnegative.'') and
$f|_{\del P}$ is the restriction of $f$ to the boundary $\del P$ of
$P$. Since the standard simplex has the structure of a compact
oriented manifold with corners, the singular chain complex of $X$
admits a natural morphism to the geometric chain complex $\Cg_\bullet
(X)$. This morphism is not a homotopy equivalence, in general, even
when $X$ is a point, as there are closed manifolds $P$ giving
nontrivial classes of oriented cobordisms. However, any geometric
chain $(P,f)$ with a boundary constraint $f(\del P) \subset A$ for
some $A \subset X$ produces a relative homology class $f_*[P] \in
H_\bullet (X, A; \Q)$, where $[P] \in H_\bullet (P, \del P; \Q)$ is
the relative fundamental class (sometimes called the ``fundamental
chain'') of $P$.

We will also need to consider geometric chains with local
coefficients. If $\F$ is a locally constant sheaf of $\Q$-vector
spaces on $X$, a \emph{geometric chain with coefficients in} $\F$ will
be a (finite) formal sum $c = \sum_i (P_i, f_i, c_i)$, where $f_i$'s
are continuous maps from compact connected oriented orbifold $P_i$'s
with corners to $X$ and $c_i$'s are global sections: $c_i \in \Gamma
(P_i; f_i^* \F)$. The differential is defined as $d c := \sum_i (\del
P_i, f_i|_{\del P_i}, c_i|_{\del P_i})$. We will use $\Cg_\bullet (X;
\F)$ to denote this complex. Since we will only consider geometric
chains, we will take the liberty to call them chains.

If $M$ is a compact connected oriented orbifold with corners, then its
\emph{fundamental chain} $[M]$ is by definition the identity map $\id:
M \to M$, understood as a geometric chain $(M, \id) \in \Cg_d (M; \Q)
= \Cg_0 (M; \Q[d])$, where $d = \dim M$ and $\Q[d]$ is the constant
sheaf $\Q$ shifted by $d$ in degree, regarded as a graded local system
concentrated in degree $-d$. If $M$ is not necessarily oriented and
$p: M^* \to M$ is the orientation cover, then we define the
\emph{fundamental chain} $[M] \in \Cg_0 (M; \Q^\epsilon)$ of $M$ to be
$(M^*, p, \frac{\ori}{2})$, where $\Q^\epsilon = \Q \times_{\mu_2}
M^*[d]$ is the \emph{orientation local system} (in particular, a
locally constant sheaf of rational graded vector spaces of rank one,
concentrated in degree $-d$) on $M$, with $M^*$ thought of as a
principle bundle over the multiplicative group $\mu_2 = \{\pm 1\}$ of
changes of orientation and $\ori \in \Gamma (M^*; p^* \Q^\epsilon)$
being the canonical orientation on $M^*$. If $M = M'/G$, where $M'$ is
an oriented compact connected orbifold with corners and $G$ a finite
group acting on $M$, then the fundamental chain of $M$ may be obtained
from the natural projection $\pi: M' \to M$ as $(M', \pi,
\frac{\ori}{\abs{G}}) \in \Cg_0 (M; \Q^\epsilon)$, where $\Q^\epsilon$
is the orientation local system of $M$. Note that a \emph{geometric
  chain with coefficients in the orientation local system} on an
orbifold $M$ may be understood as a linear combination of geometric
chains $f: P \to M$ with a (continuous) choice of \emph{local
  orientation on $M$ along} $P$.


\section{The Open-Closed Moduli Space}
\label{OCM}

In this section, we will describe a certain moduli space $\MUg$, which
will be used to set up and solve the Quantum Master Equation (QME).

That space $\MUg$ will be closely related to the moduli space $\MM$,
introduced by Melissa Liu \cite{liu}, of stable bordered Riemann
surfaces of type $(g,b)$ with $(n,\v{m})$ punctures (marked points).
The space $\MM$ generalizes the \emph{Deligne-Mumford moduli space}
$\overline{\M}_{g,n}$, parameterizing isomorphism classes of stable
algebraic curves of genus $g$ with $n$ punctures, to the open-closed
case, i.e., $\MM$ parameterizes isomorphism classes of stable bordered
Riemann surfaces with $b$ boundary components, $n$ punctures in the
interior, and $m_i$ punctures on the $i$th boundary component, if
$\v{m} = (m_1, \dots, m_b)$, see more details in the next paragraph.
All the boundary components are considered labeled by numbers 1
through $b$, and all the punctures must be distinct and labeled 1
through $n$ for the interior punctures and 1 through $m_i$ for the
punctures on the $i$th boundary component, $i = 1, \dots, b$.  We
require that labels on the boundary punctures must be placed in a way
compatible with the induced orientation on the boundary, i.e., the
cyclic, $\mod m_i$ order of punctures on the $i$th boundary component
must increase when moving along the boundary \emph{counterclockwise},
that is to say, keeping the surface on the right-hand side. In the
ideal world conveniently provided by string theory, interior and
boundary punctures allegedly correspond to evolving closed and open
strings, respectively, whose world sheet is the Riemann surface.

A \emph{bordered Riemann surface} here means a complex curve with real
boundary, i.e., a compact, connected, Hausdorff topological space,
locally modeled on the upper half-plane $H = \{ z \in \C \; | \; \im z
\ge 0\}$ using analytic maps.  A \emph{prestable bordered Riemann
  surface} is a bordered Riemann surface with at most a finite number
of \emph{singularities of nodal type} at points other than the
punctures. The \emph{allowed types of nodes} are denoted X, E, and H,
where X means an interior node (locally isomorphic to a neighborhood
of 0 on $\{ xy = 0 \} \subset \C^2$), E a boundary node, when the
whole boundary component is shrunk to a point (locally modeled on a
neighborhood of 0 on $\{x^2 + y^2 = 0 \subset \C^2\}/ \sigma$, where
$\sigma (x,y) = (\bar x, \bar y)$ is the complex conjugation), and H a
boundary node at which a boundary component intersects itself or
another boundary component (locally modeled on a neighborhood of 0 on
$\{x^2 - y^2 = 0 \subset \C^2\}/ \sigma$).

A prestable bordered Riemann surface is \emph{stable}, if its
automorphism group is discrete. Here an automorphism must map the
boundaries to the boundaries and the punctures to the punctures,
respecting the labels. The stability condition is equivalent to the
condition that the Euler characteristic (in a certain generalized
sense, see below) of each component of the surface obtained by
removing all the punctures is negative. Note that the \emph{Euler
characteristic} is by definition one half of the Euler characteristic
of the double. Thus, for a nondegenerate bordered surface $\Sigma$,
its Euler characteristic is given by
\begin{equation}
\label{euler}
 \chi (\Sigma) = 2- 2g - b - n - m/2.
\end{equation}
The stability condition thereby excludes a finite number of types
$(g,b; n,\v{m})$, namely $g=b=0$ with $n \le 2$; $g=1$, $b = 0$ with
$n=0$; $g=0$, $b=1$ with $n \le 1$, $m=0$ or $n=0$, $m \le 2$; and
$g=0$, $b=2$ with $m=n=0$. The spaces $\MM$ have been thoroughly
studied by M.~Liu in \cite{liu}. They are compact, Hausdorff
topological spaces with the structure of a smooth orbifold with
corners of dimension $6g-6+2n+3b+m$, where $m = \sum_{i=1}^b m_i$ is
the total number of boundary punctures.

Our space $\MUg$ is the moduli space of isomorphism classes of stable
bordered Riemann surfaces of type $(g,b)$ with $(n,\v{m})$ punctures
and certain extra data, namely, decorations by a real tangent
direction, i.e., a ray, in the complex tensor product of the tangent
spaces on each side of each interior node. The space $\MUg$ can be
obtained by performing real blowups along the divisors of $\MM$
corresponding to the interior nodes, as in \cite{ksv}. The dimension
of $\MUg$ is the same as that of $\MM$: $\dim \MUg = \dim \MM =
6g-6+2n+3b + m$.

We will concentrate on the moduli space
\[
\MUgm/\S = \left. \left(
  \coprod_{\v{m}: \sum m_i = m} \MUg / \Z_{m_1} \times \dots \times \Z_{m_b}
\right) \right/ \S_b \times \S_n
\]
of stable bordered Riemann surfaces as above with \emph{unlabeled}
boundary components and punctures, that is, the quotient of the
disjoint union $\displaystyle{\MUgm = \coprod_{\v{m} : \sum m_i = m}
  \MUg}$ of moduli spaces with labeled boundaries and punctures by an
appropriate action of the permutation group $\displaystyle{\S = \left(
    \prod_{\v{m}: \sum m_i = m} \Z_{m_1} \times \dots \times \Z_{m_b}
  \right) \rtimes \S_b \times \S_n}$. To set up the QME, we will work
with geometric chains of this moduli space with twisted coefficients,
i.e., a one-dimensional local system $\Q^\epsilon$ obtained from a
certain sign representation $\rho: \S \to
\End L = \Q^*$ of the permutation group $\S$ in a one-dimensional
rational graded vector space $L$ concentrated in degree $-d := - \dim
\MUgm$. This representation is defined as follows: $\rho$ is a trivial
representation of $\S_n$;
\[
 \rho (\zeta_i) = (-1)^{m_i - 1}
\]
for the generator $\zeta_i (p) = p+1 \mod m_i$ of the group $\Z_{m_i}$
of cyclic permutations of the punctures on the $i$th boundary
component (where $m_i \ge 1$); and
\[
 \rho(\tau_{ij}) = (-1)^{(m_i-1)(m_j-1)}
\]
for the transposition $\tau_{ij} \in \S_b$ interchanging (the labels
of) the $i$th and $j$th boundary components. Then $\Q^\epsilon$ is the
locally constant sheaf $\MUgm \times_\S L$ over $\MUgm /\S$. Note that
an ordering of the boundary components and an ordering of the boundary
punctures on each boundary component compatible with the cyclic
ordering thereof on a given Riemann surface determines a section of
the local system over the point in the moduli space $\MUgm/\S$
corresponding to that Riemann surface. A change of these orderings
will change this section by a sign factor as defined by the
representation $\rho$.

\section{Orientation}

Here we claim that the space $\MUg$ is an orientable orbifold with
corners.  Recall that an orientation on an orbifold with corners is
the choice of a nowhere vanishing section of the orbifold determinant
tangent (or, equivalently, cotangent) bundle up to a positive real
function factor. An orientation on $\MUg$ may be defined similarly to
an orientation on the Deligne-Mumford-Liu space $\MM$, see
\cite[Theorem 4.14]{liu}, as follows.

The interior (i.e., pre-compactification) moduli space
$\Mgm{1,\dots,1}$ is an orbifold with a natural complex structure, see
\cite{ivash-shevch}, and thereby has a natural orientation. The reason
is that this moduli space is isomorphic to the moduli space of complex
algebraic curves with $n$ labeled punctures with a \emph{holomorphic}
involution of a certain topological type (namely, with $b$ invariant
closed curves), which is naturally complex as a moduli space of
complex objects. We will orient the spaces $\Mg$ for $\v{m} = (m_1,
\dots, m_b)$, $m_i \ge 0$, inductively, by lifting or projecting
orientation along the following fiber bundles:
\begin{equation}
\label{forget}
\Mgm{m_1,\dots,m_i+1, \dots, m_b} \to \Mgm{m_1,\dots,m_i, \dots, m_b},
\end{equation}
which have naturally oriented fibers, identified with the open arc of
the $i$th boundary component $B_i$ of the Riemann surface between the
first and the $m_i$th punctures. As usual, we use the counterclockwise
orientation on the boundary components under which the surface is
always on the right-hand side. We say that the orientations on
$\Mgm{m_1,\dots,m_i+1, \dots, m_b}$ and $\Mgm{m_1,\dots,m_i, \dots,
  m_b}$ agree, if the frame of tangent vectors to
$\Mgm{m_1,\dots,m_i+1, \dots, m_b}$ obtained by appending a
counterclockwise tangent vector to $B_i$ to a positively oriented
frame of tangent vectors to $\Mgm{m_1,\dots,m_i, \dots, m_b}$ is
positively oriented. To endow all the spaces $\Mg$, $m_i \ge 0$ for
all $i$, with orientation, start with the orientation on
$\Mgm{1,\dots,1}$ coming from the complex structure and order the
total of $m \ge 0$ open punctures in a linear fashion: 1 through $m_1$
(in the order in which they are labeled) on the first boundary
component, then $m_1+1$ through $m_1+m_2$ on the second boundary
component, and so on. Then use the fiber bundles \eqref{forget}, one
by one, in the resulting order, to induce orientation on all spaces
$\Mg$ and thereby on its compactification $\MUg$ to an orbifold with
corners.

The choice of orientation determines uniquely a fundamental chain
$[\MUg] \linebreak[0] \in \linebreak[1] \Cg_{6g-6+2n+3b+m} (\MUg; \Q)$
of the orbifold $\MUg$ with corners. We are rather interested in the
unlabeled moduli space $\MUgm/\S$, though.

\begin{prop}
  The local system $\Q^\epsilon$ is the orientation sheaf for the
  orbifold $\MUgm/\S$ with corners.
\end{prop}

\begin{proof}
  The (real) orientation bundle of an orbifold $M$ with corners is by
  definition given by the determinant tangent bundle of $M$, regarded
  as an orbifold real vector bundle. Here by the determinant tangent
  bundle we mean the top graded symmetric power $\det TM := S^d
  (TM[1])$ of the tangent bundle $TM$ placed in degree $-1$, where $d
  = \dim M$.  The real orientation bundle is induced from a unique (up
  to isomorphism) locally constant sheaf of graded $\Q$-vector spaces
  by extension of scalars. Thus, it is enough to talk about the real
  orientation bundle in the proof.

  The fact that $\MUg$ is orientable means the determinant bundle
  $\det \linebreak[0] T \MUg$ of the tangent bundle of its interior
  part $\MUg$ is trivial. The orientation bundle over $\MUgm/\S$ will
  then be $\det T (\MUgm/\S)$. Looking at how orientation was defined
  on $\MUg$, observe that on the union $\coprod_{\v{m} : \sum m_i = m}
  \det T \MUg$ of bundles, the permutation group $\S_n$ of the set of
  interior punctures acts trivially, while transposition of the $i$th
  and the $j$th boundary components acts by $(-1)^{(m_i-1)(m_j-1)}$,
  and the basic cyclic permutation $\zeta_i$ acts by $(-1)^{m_i-1}$,
  if $m_i \ge 1$.
\end{proof}

Thus, the \emph{fundamental chain} $[\MUgm/\S] \in \Cg_0 (\MUgm/\S;
\Q^\epsilon)$ is well defined, see Section~\ref{chains}.

\section{BV Structure for the Open-Closed Moduli Spaces}
\label{BV}

A \emph{differential graded $($dg$)$ Batalin-Vilkovisky $($BV$)$
  algebra} structure on a complex $V$ (i.e., a vector space with a
differential $d: V \to V$, $d^2 = 0$, of degree 1) consists of a
(graded) commutative \emph{dot product} $V \otimes V \to V$, $a
\otimes b \mapsto ab$, and a \emph{BV operator} $\Delta: V \to V$
which is a second-order differential of degree 1. Both operations must
be compatible with the differential $d$, which must be a (graded)
derivation of the dot product (assuming annihilation of constants) and
(graded) commute with the BV operator: $[d,\Delta] := d \Delta +
\Delta d = 0$. Secondary to this basic structure is an
\emph{antibracket}
\[
\{ a, b\} := (-1)^{\abs{a}} \Delta (ab) - (-1)^{\abs{a}} \Delta (a) b
- a \Delta b,
\]
which turns out to be a Lie bracket of degree 1 on which both
differentials $d$ and $\Delta$ act as derivations. In fact, the
derivation property of the BV operator with respect to the antibracket
is equivalent to its second-order derivation property with respect to
the dot product.

Consider the space
\[
U = \bigoplus_{g,b,m, n} \Cg_\bullet (\MUgm/\S; \Q^\epsilon)
\]
of geometric chains.
Since the coefficient system carries a degree shift by the dimension
of $\MUgm/\S$, the space $U$ carries a natural grading by the negative
codimension of the chain in $\MUgm/\S$. We will take the opposite
grading on $U$, i.e., the \emph{grading by the codimension} of the
chain in $\MUgm/\S$. The differential $d$ of geometric chains will then
have degree 1 and make $U$ into a complex of rational vector spaces.

The space $V$ on which we will introduce a dg BV algebra structure
will be defined as follows:
\[
V := \bigoplus_{b, m, n} \Cg_\bullet (\MUm/\S; \Q^\epsilon),
\]
where $\MUm/\S$ is the moduli space of stable bordered Riemann
surfaces with $b$ boundary components, $n$ interior punctures, and $m$
boundary punctures, just like the Riemann surfaces in $\MUgm/\S$, but
in general having multiple connected components of various genera. The
grading on $V$ is given by codimension in the corresponding connected
component of $\MUm/\S$, and the dot product is induced by disjoint
union of Riemann surfaces. We formally add a copy of the ground field
$\Q$ to $V$, and the unit element $1 \in \Q \subset V$ might be
interpreted as the fundamental chain of the one-point moduli space
comprised by the empty Riemann surface. The algebra $V$ is also
isomorphic to the graded symmetric algebra $S(U)$.  We will be using
this observation when discussing the $L_\infty$ structure on $U$
later.

To define a dg BV algebra structure on $V$, it remains to define a BV
operator $\Delta: V \to V$ satisfying required properties. It will
consist of three components:
\[
\Delta = \Delta_c + \Delta_o + \delco,
\]
each of degree 1, square zero, and (graded) commuting with each other.

The operator $\Delta_c$ is induced on chains by twist-attaching at
each pair of interior punctures.  To achieve this, we define a bundle
$ST \MUm/\S$ over $\MUm/\S$ of triples $(\Sigma, P, r)$, where $\Sigma
\in \MUm/\S$, $P$ is a choice of an unordered pair of interior
punctures on $\Sigma$, and $r$ is one of the $S^1$ ways of attaching
them (i.e., a real ray in the tensor product (over $\C$) of the
tangent spaces of $\Sigma$ at these two punctures). So the fiber is
homeomorphic to $\coprod^{n(n-1)/2} S^1$.  Then we have a diagram
\begin{equation}
\label{closedpp}
\MUm/\S \xleftarrow{\pi} ST \MUm/\S \xrightarrow{a_c}
\MUmnv{n-2}{m}{b}/\S,
\end{equation}
where $\pi$ is the bundle projection map and $a_c$ is obtained by
attaching the two chosen punctures $P$ on $\Sigma$ and decorating the
resulting node with the chosen real ray $r$.  (One can view this
diagram as a morphism realizing twist-attaching in the category of
correspondences.)  Then \emph{twist-attaching} for chains is defined
as the corresponding ``push-pull,'' giving us the ``closed'' part
$\Delta_c$ of the BV operator:
\begin{equation*}
\Delta_c := (a_c)_* \pi^!: \Cg_\bullet (\MUm/\S; \Q^\epsilon)
\to \Cg_{\bullet+1} (\MUmnv{n-2}{m}{b}/\S; \Q^\epsilon),
\end{equation*}
Here the pullback $\pi^!$ for geometric chains is simply the geometric
pre-image. More precisely, to define the pullback of a geometric chain
$(P,f,c)$, we take the pullback $f^* ST$ of the fiber bundle $ST$
along $f$ and the chain $(f^* ST, \tilde f)$, where $f^* ST$ is the
total space and $\tilde f: f^* ST \to ST$ is the pullback of $f$ (if
$f^*ST$ is disconnected we regard $\tilde f$ as a sum of maps).  To
define what $\pi^!$ does to a section $c$ of $\Q^\epsilon$, lift the
diagram \eqref{closedpp} to

\begin{equation*}
\MUm \xleftarrow{\pi} ST \MUm \xrightarrow{\Delta_c}
\MUmnv{n-2}{m}{b},
\end{equation*}

\noindent defined before taking the quotient by the symmetric groups.
Here, $ST \MUm$ is the bundle whose fiber over $\Sigma \in \MUm$
consists of all the $n(n-1)/2$ possible choices of unordered pairs
$\{i,j\}$ of labeled punctures along with the $S^1$ ways of attaching
them.  The fiber of $\pi$ (isomorphic to $n(n-1)/2$ copies of $S^1$)
has a natural orientation coming from the counterclockwise orientation
on the tensor product (over $\C$) of the tangent spaces at the
punctures $i$ and $j$.  The orientation on the total space $ST \MUm$
is then defined (locally) as the orientation on the base $\MUm$ times
the orientation of the fiber of $\pi$, and the orientation sheaf on
$ST \MUm/\S$ determines a local system.  Now recall that a section $c$
of $\Q^\epsilon$ is a rational number $c'$ multiplied by the
orientation $\ori$ of the moduli space $\MUm$. Since an orientation of
$\MUm$ determines an orientation on $ST \MUm$, as we have just
described, we use that orientation, multiplied by the same number
$c'$, to get a section of the local system on $ST \MUm/\S$.

The operator $\Delta_o$ is induced on geometric chains by attaching at
each pair of boundary punctures. To describe this procedure precisely,
we form the bundle $B' \MUm/\S$ where the fiber over a point $\Sigma
\in \MUm/\S$ consists of all possible choices of pairs of punctures on
$\Sigma$ which both lie on the same boundary component.  Similarly,
form the bundle $B'' \MUm/\S$ whose fibers are all possible pairs
of punctures lying on different boundary components. Then we get the
following diagrams:
\[
\MUm/\S \xleftarrow{\pi'} B' \MUm/\S \xrightarrow{a'_o}
\MUmnv{n}{m-2}{b+1}/\S
\]
and
\[
\MUm/\S \xleftarrow{\pi''} B'' \MUm/\S \xrightarrow{a_o''}
\MUmnv{n}{m-2}{b-1}/\S,
\]
where $a'_o$ and $a_o''$ are the obvious attaching maps.

We perform the push-pull again to obtain chain maps $\Delta'_o$ and
$\Delta_o''$.  The bundles before quotienting, $B' \MUm$ and $B''
\MUm$, are just direct products of $\MUm$ with finite discrete sets
and are thus orientable, and we define the pullback (in this case also
known as the ``transfer homomorphism'') of geometric chains as in the
closed case.  The pushforward of sections of the local system is
defined in the next paragraph.  We will then define the corresponding
component of the BV operator as
\[
\Delta_o := \Delta'_o + \Delta_o''.
\]

\begin{figure}
\includegraphics[width=7cm]{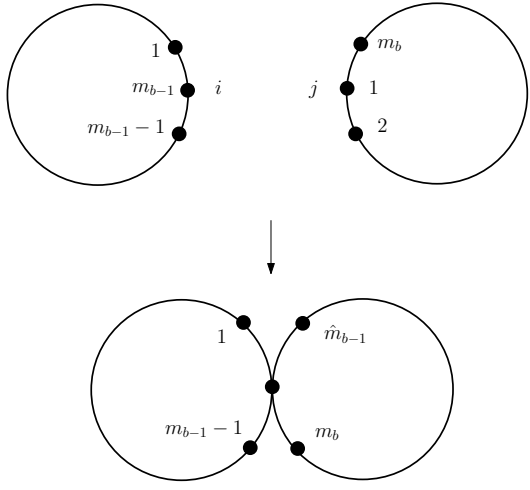}
\caption{Attaching two punctures on different boundary components.
The resulting new boundary component has $\hat{m}_{b-1} = m_{b-1} + m_b - 2$
punctures.}
\label{zwiebach3a}
\end{figure}

Now let us define the pushforward of sections of the local system via
$a'_o$ and $a_o''$.  Recall that an ordering of the boundary
components and a cyclic ordering of the punctures on each boundary,
for a given surface $\Sigma \in \MUm/\S$, gives a section of the local
system over the point $\Sigma$. The same can be said for $(\Sigma, P)
\in B' \MUm/\S$ or $B'' \MUm/\S$, where $P$ is the choice of a pair of
boundary punctures on $\Sigma$.  Thus to define the pushforward it
suffices to explain how attaching acts on the labeling of boundaries
and boundary punctures.  If the punctures $i$ and $j$ in the pair $P$
lie on different boundary components, see Figure~\ref{zwiebach3a},
first change the ordering of the boundary components and boundary
punctures in a way that the puncture $i$ is the last puncture on the
$b-1$st boundary component and the puncture $j$ is the first puncture
on the $b$th boundary component. Then, after the punctures are
attached to form a single boundary component, order the boundary
components so that this new one goes last (i.e., becomes number
$b-1$), with the same ordering of the old boundary components. Order
the punctures on the new boundary component by placing the punctures
coming from the old $b-1$st boundary component first, preserving their
order, followed by the punctures coming from the old $b$th boundary
component, in their old order. Keep the old ordering of punctures on
the boundary components not affected by attaching.

\begin{figure}
\includegraphics[width=10cm]{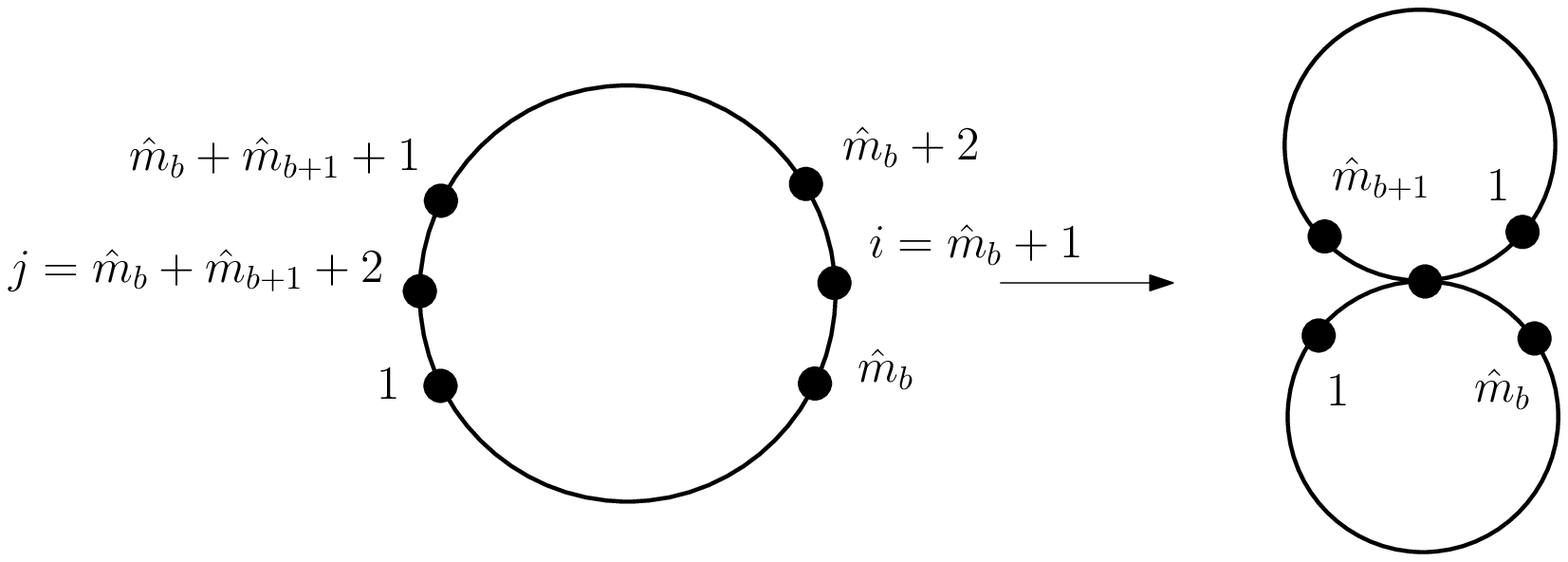}
\caption{Attaching two punctures on the same boundary component with
$m_b = \hat{m}_b + \hat{m}_{b+1} + 2$ punctures.}
\label{pinch}
\end{figure}

In the case when the punctures $i$ and $j$ happen to be on the same
boundary component, see Figure~\ref{pinch}, change the ordering of the
boundary components so that this component goes last and the puncture
$j$ is the last puncture on that boundary component. Out of the two
boundary components obtained by the pinching, order the one following
$j$ counterclockwise first, the other new boundary component next,
preceded by the old boundary components in the old order. Order the
punctures on the two new boundary components, declaring the puncture
going after the double point on the first new boundary component to be
first, followed by the other punctures in the counterclockwise order,
and the puncture going after the double point on the second boundary
component in the counterclockwise manner first on that boundary
component. Again, keep the old ordering of punctures on the boundary
components not affected by attaching.

Finally, define $\delco$ as follows. Let $B \MUm/\S$ be the bundle
over $\MUm/\S$ whose fiber over $\Sigma$ consists of all choices of an
interior puncture on $\Sigma$. Then we get:
\[
\MUm/\S \xleftarrow{\pi} B \MUm/\S \xrightarrow{a_{\operatorname{co}}}
\MUmnv{n-1}{m}{b+1}/\S
\]
where $a_{\operatorname{co}}$ is the map induced by declaring the
chosen interior puncture to be a degenerate boundary component with no
``open'' punctures on it.  Then the chain map $\delco$ is defined as
the corresponding push-pull where the pushforward of the section of
our local system is defined by putting this boundary component after
the other boundary components in the ordering.

\begin{thm}
  The operator $\Delta = \Delta_c + \Delta_o + \delco$ is a
  graded second-order differential on the dg graded commutative
  algebra $V$ and thereby defines the structure of a dg BV algebra on
  $V$.
\end{thm}

\begin{proof}
  The fact that $\Delta (1) = 0$ follows tautologically from the
  definition of the three components of $\Delta$. It will be enough to
  check the following identities:
\begin{gather*}
  [\Delta_c, d] = [\Delta_o, d] = [\delco, d] = 0,\\
  \Delta_c^2 = \Delta_o^2 = \delco^2  = 0,\\
  [\Delta_c, \Delta_o] = [\Delta_c, \delco]
  = [\Delta_o, \delco]  = 0,\\
  \Delta_c \text{ and } \Delta_o  \text{ are second-order differential operators},\\
  \delco \text{ is a derivation}.
\end{gather*}
The operator $\Delta_c$ and the differential $d$ commute, because the
copy of $S^1$ acquired by twist-attaching is a closed manifold. The
commutation of $d$ with the other components of $\Delta$ is more
obvious.

The fact that $\Delta_c$ is a differential, i.e., $\Delta_c^2 = 0$,
comes from our definition of orientation: we place the extra component
$S^1$ last in orientation, so that each term in $\Delta_c^2 C$ will
have two extra twists $S^1$, as compared to the original chain $C$,
and will be canceled by another term in $\Delta_c^2 C$, in which the
two twists come in the opposite order. The property $\Delta_o^2 = 0$
is also true because of our choice of orientation. Each term in
$\Delta_o^2 C$ is obtained by attaching one pair of boundary punctures
together and then another pair. This term will be canceled by the term
in which those two pairs of punctures are attached in the opposite
order. It is a straightforward calculation to see that the signs
coming from the our local system work out to cancel those pairs of
terms in $\Delta_o^2 C$.

A more conceptual explanation of the same phenomenon may be done using
the interpretation of the local system as the orientation sheaf of our
orbifold. Note that the choice of orientation under attaching a pair
of open punctures is performed in a way that we remove one factor in
the top wedge power of the tangent bundle to the moduli space $\MUm$,
leaving the other factors intact. The corresponding pairs of terms in
$\Delta_o^2 C$ in which the same two pairs of punctures are attached
in the opposite order will cancel each other, because the orders in
which the corresponding factors are removed from the wedge product
will be opposite. The same argument applies to showing $\delco^2 = 0$
and, in fact, the graded commutation
\[ 
  [\Delta_c, \Delta_o] = [\Delta_c, \delco]
  = [\Delta_o, \delco]  = 0.
\]

The fact that $\delco$ is a graded derivation of the dot product,
\[
\delco (a \cdot b) = \delco (a) \cdot b + (-1)^\abs{a} a\cdot
\delco(b),
\]
is obvious: transformation of an interior puncture into a degenerate
boundary component on a disjoint union of two Riemann surfaces happens
on either one surface or the other.

The fact that $\Delta_c$ is a second-order derivation is equivalent to
the following statement. Define a bracket
\[
\{ a, b\}_c := (-1)^{\abs{a}} \Delta_c (ab) - (-1)^{\abs{a}} \Delta_c (a) b
- a \Delta_c b.
\]
Then this bracket is a graded derivation in each (or equivalently,
one) of its variables, that is
\begin{equation}
\label{leibniz}
\{a, bc\}_c = \{a,b\}_c c + (-1)^{(\abs{a}+1)\abs{b}} b \{a,c\}_c.
\end{equation}

What is clear from the definition, the geometric meaning of the
bracket $\{a,b\}_c$ is $(-1)^{\abs{a}}$ multiplied by the alternating
sum of twist-attachments over all pairs of closed punctures for the
chains $a$ and $b$, respectively. Given three geometric chains $a$,
$b$, and $c$, Equation \eqref{leibniz} is obvious, as it just says
that twist-attaching of closed punctures in $a$ with those in the
disjoint union of $b$ and $c$ breaks into twist-attaching with
punctures in $b$ and $c$ and then taking the disjoint union. The signs
come out right, because of our definition of orientation under
twist-attaching and disjoint union.

The same argument applies to $\Delta_o$. Consider a bracket
\[
\{ a, b\}_o := (-1)^{\abs{a}} \Delta_o (ab) - (-1)^{\abs{a}} \Delta_o (a) b
- a \Delta_o b
\]
and show it satisfies the derivation property
\[
\{a, bc\}_o = \{a,b\}_o c + (-1)^{(\abs{a}+1)\abs{b}} b \{a,c\}_o.
\]
\end{proof}

\begin{rem}
  Note that the part $\Delta'_o$ of $\Delta_o$ corresponding to
  attaching punctures lying on the same boundary component is actually
  a derivation, and therefore the open part of the antibracket comes
  only from $\Delta''_o$ corresponding to attaching punctures lying on
  different boundary components:
\[
\{a,b\}_o' = 0, \qquad \{ a, b\}_o = \{ a, b\}''_o.
\]
\end{rem}

\section{Quantum Master Equation}
\label{section:QME}

In any dg BV algebra, it makes sense to set up the \emph{Quantum
Master Equation} (\emph{QME}):
\begin{equation}
\label{QME1}
dS + \Delta S + \frac{1}{2}\{S,S\} = 0
\end{equation}
for $S \in V$ of degree zero, in which case all the terms will be in
the same degree. If we allow the formal power series
\[
e^S = 1 + S + S^2/2! + S^3/3! + \dots,
\]
then, using the second-order differential operator property of
$\Delta$, the QME may be written equivalently as follows:
\begin{equation}
\label{QME2}
(d + \Delta) e^S = 0.
\end{equation}

In our context, we will need an extra formal variable $\lambda$ of
degree 0, called the \emph{string coupling constant}, and set up the
QME in the space:
\[
V[[\lambda]] : = \{ \sum_{n=0}^\infty v_n \lambda^n \; | \; v_n \in V\},
\]
which inherits a dg BV algebra structure from $V$ by linearity in
$\lambda$. More generally, we will consider the following modification
of the QME in the dg BV algebra $V[[\lambda,\sqrt{\hslash}]]$, where
$\hslash$ is another formal variable of degree 0, called the
\emph{Planck constant}:
\begin{equation}
\label{MQME1}
dS + \hslash \Delta S + \frac{1}{2} \{ S, S\} = 0,
\end{equation}
or, equivalently, in $V[[\lambda,\sqrt{\hslash},
(\sqrt{\hslash})^{-1}]]$
\begin{equation}
\label{MQME2}
(d + \hslash \Delta) e^{S/\hslash} = 0.
\end{equation}
These last equations turn into the QMEs \eqref{QME1} and \eqref{QME2},
respectively, under the specification $\hslash = 1$, so that if $S
(\lambda, \hslash) \in V [[\lambda, \sqrt{\hslash}]]$ is a solution of
\eqref{MQME1} and $S(\lambda, 1)$ makes sense, then it will
automatically satisfy \eqref{QME1}. We will deal with the more general
QME of the form \eqref{MQME1} or \eqref{MQME2} in this paper. However,
we will modify the BV operator on $V[[\lambda, \sqrt{\hslash}]]$ in
our case, when $V$ is the space of geometric chains of the open-closed
moduli space, in the following way:
\[
\Delta := \Delta_c + \Delta_o + \sqrt{\hslash} \delco.
\]
Note that since $\delco$ is a first-order differential operator, the
new $\Delta$ will still be a BV operator on the algebra $V[[\lambda,
\sqrt{\hslash}]]$ over the ring $\Q[[\lambda, \sqrt{\hslash}]]$, and
the modification does not change the antibracket. Note also the change
does not affect $\Delta$ after we make the evaluation $\hslash =1$.

Take
\begin{equation}
\label{solution}
S : = \sum_{g,b, m, n} S_{g,b}^{n,m} \lambda^{-2\chi} \hslash^{p-\chi},
\end{equation}
where $S_{g,b}^{n,m} := [\MUgm/\S] \in \Cg_0 (\MUgm/\S; \Q^\epsilon)$
denotes the fundamental chain, $\chi$ is the Euler characteristic
\eqref{euler} of the bordered Riemann surface representing a point in
$\MUgm$, the summation runs over the indices corresponding to stable
moduli spaces, i.e., $\chi = 2 - 2g - b - n - m/2 < 0$, and
\[
p = 1- (m+n)/2.
\]

\begin{thm}
S is a solution of the QME, i.e.,
\[
(d + \hslash \Delta) e^{S/\hslash} = 0.
\]
\end{thm}

\begin{proof}
  We will prove the equation in the equivalent form \eqref{MQME1}.
  This equation in fact encodes the structure of the boundary $\del
  \MUgm/\S$ as the result of the operation on lower-dimensional moduli
  spaces induced by attaching open-closed Riemann surfaces at
  punctures and regarding interior punctures as degenerated boundary
  components. We will start with representing the boundary of a single
  moduli space $\MUgm/\S$ in this way and then pass to the generating
  series~\eqref{solution} to get the QME.

  The boundary of the open-closed moduli space $\MUgm$ is given by the
  locus of Riemann surfaces with at least one node. There are three
  types of nodes, X, E, and H, see Section~\ref{OCM}. Type X and E
  nodes may also be subdivided into two types depending on whether the
  node separates two irreducible components of the Riemann surface.
  Likewise, nonseparating type E nodes split into two types, depending
  on whether they are obtained by pinching together two punctures on
  the same or different boundary components. This implies that the
  boundary of the moduli space $\MUgm/\S$ with unlabeled punctures is
  given exactly as follows, see also Figures \ref{formula1} and
  \ref{formula2}:
\begin{figure}
\includegraphics[width=10cm]{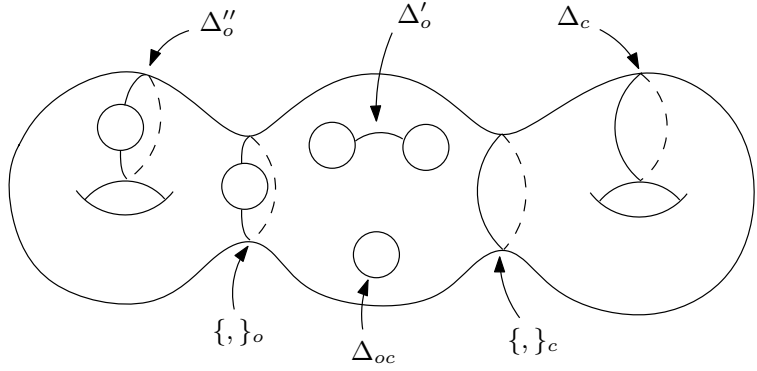}
\caption{Different types of degeneration, resulting from contraction of
curves on the surface.}
\label{formula1}
\end{figure}
\begin{figure}
\includegraphics[width=12cm]{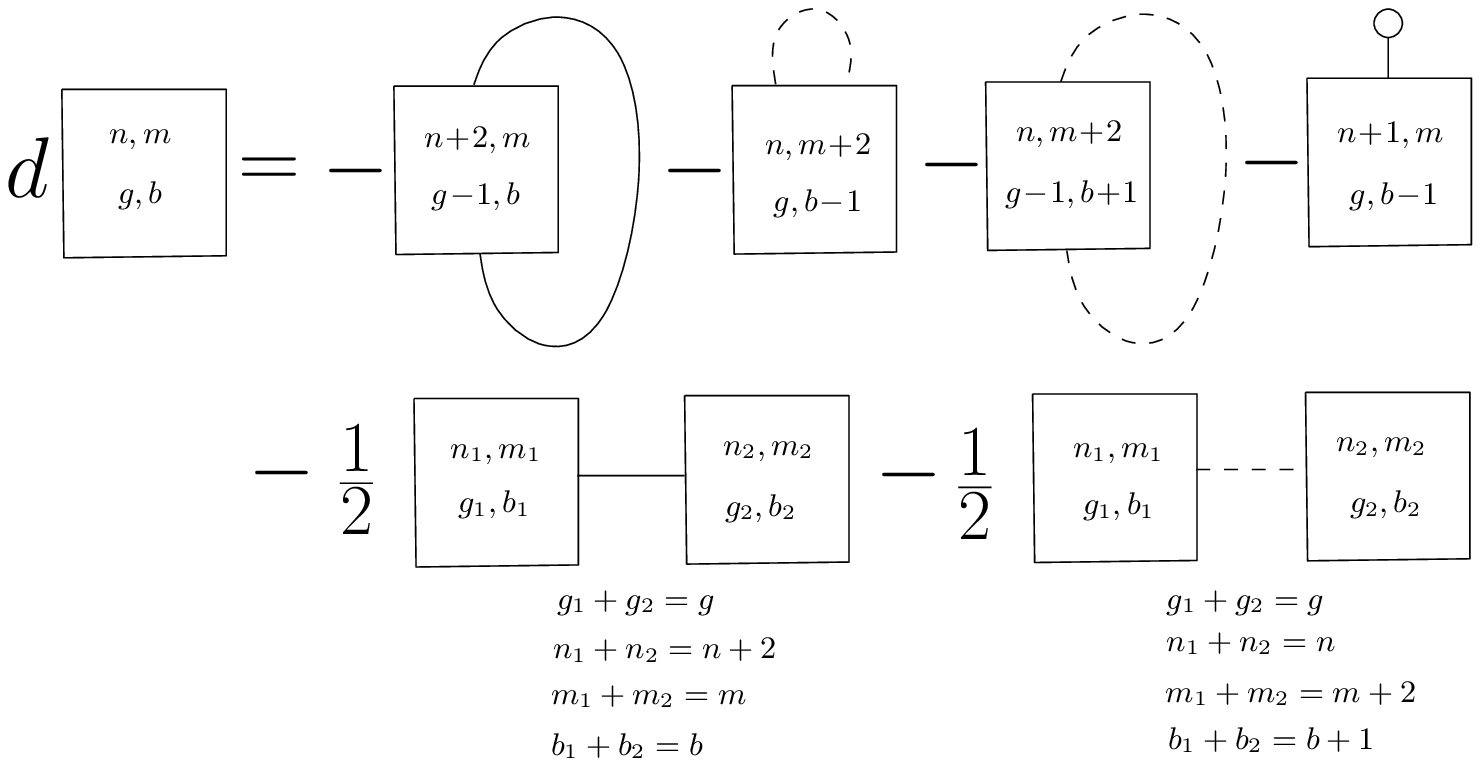}
\caption{The boundary of the moduli space: the solid connector means
twist-attaching at interior punctures, the dotted connector means
attaching at boundary punctures, the circle on a pole means turning an
interior puncture into a degenerated boundary component.}
\label{formula2}
\end{figure}
\begin{multline}
\label{premaster}
d S_{g,b}^{n,m} = -\Delta_c S^{n+2,m}_{g-1,b} - \Delta'_o
S^{n,m+2}_{g,b-1} - \Delta''_o S^{n,m+2}_{g-1,b+1} - \delco
S^{n+1,m}_{g,b-1}
\\
- \frac{1}{2} \sum_{\substack{
g_1 + g_2 = g, \; b_1 + b_2 = b,\\
n_1 + n_2 = n+2, \; m_1 + m_2 = m}}
\{S^{n_1,m_1}_{g_1,b_1} , S^{n_2,m_2}_{g_2,b_2}\}_c
\\
- \frac{1}{2} \sum_{\substack{
g_1 + g_2 = g, \; b_1 + b_2 = b+1,\\
n_1 + n_2 = n, \; m_1 + m_2 = m+2}}
\{S^{n_1,m_1}_{g_1,b_1} , S^{n_2,m_2}_{g_2,b_2}\}_o.
\end{multline}
The $-1$ and $-1/2$ factors can be explained in the following way. The
negative signs on the right-hand side are due to our choice of
orientation. For example, the component of the boundary $\del \MUgm$
corresponding to contracting a closed real curve in the interior of
the Riemann surface contributes to $\Delta_c S^{n+2,m}_{g-1,b}$ on the
right-hand side with a negative sign, because the corresponding
angular Fenchel-Nielsen coordinate goes last in the orientation
defined for $\Delta_c$, while the boundary on the left-hand side is
oriented in the way that the last coordinate, which is the
circumference Fenchel-Nielsen coordinate, is nonnegative. Otherwise,
all the remaining coordinates may be chosen with the same order on
both sides. Since in the complex structure in Fenchel-Nielsen
coordinates, the perimeter coordinate goes just before the angular
coordinate for the same real curve, the perimeter coordinate is
actually the second-to-last in the orientation of $\MUgm$, whence the
negative sign.

Let us look at another example, the $\Delta'_o S^{n,m+2}_{g,b-1}$ term
in Equation~\eqref{premaster}. Again, locally on the moduli space, if
the last two boundary components (if not the last two, we will change
the ordering, changing orientation on the moduli space, if necessary),
number $b-1$ and number $b$, get close together and touch, the moduli
space will be given by the equation $r \ge r_0$, where $r$ is the
distance between the centers of the two circles. This $r$ is the first
real coordinate in the polar coordinate system $(r,\theta)$, where
$\theta$ is the angle of the direction from the center of the $b-1$st
circle to the center of the $b$th circle. This coordinate system has
the same orientation as the one coming from the complex coordinate
given by the position of the center of the $b$th circle near the
$b-1$st one on the Riemann surface. When the orderings of the
punctures on the circles are chosen in a way that the circles touch at
points between the first and the last punctures on each circle (we can
always do that by changing the orderings and adjusting orientation),
this ordering will be the same as the one coming from the operation
$\Delta'_o$ of attaching two punctures on the $b-1$st boundary
component of the moduli space $\MUmnv{n}{m+2}{g,b-1}$, as described by
Figure~\ref{pinch}, except that the figure refers to $\Delta'_o$ on
$\MUgm$. This all gives the same orientation as on the boundary of
$\MUgm$, but since an oriented orbifold with positively oriented
boundary must be given by the last coordinate being nonnegative and
$r$ was the second-to-last coordinate in a positively oriented
coordinate system, we get a negative sign in front of $\Delta'_o
S^{n,m+2}_{g,b-1}$ in Equation~\eqref{premaster}.

The factor of $1/2$ is due to the fact that each term
$\{S^{n_1,m_1}_{g_1,b_1}, S^{n_2,m_2}_{g_2,b_2}\}_c$ in the sum for
the closed part of the antibracket, is counted twice (in the ``orbi''
sense), the second time as $\{S^{n_2,m_2}_{g_2,b_2},
S^{n_1,m_1}_{g_1,b_1}\}_c$, even though this term is present only once
in the boundary on the left-hand side.

\begin{sloppypar}
  Now as we have checked Equation~\eqref{premaster}, let us sum up
  these equations over different $g$, $b$, $m$, and $n$ with weights
  $\lambda^{-2\chi} \hslash^{p-\chi}$ in a single equation for the
  generating series. Then use the fact that the Euler characteristic
  of an open-closed Riemann surface does not change under degeneration
  of the surface, or equivalently, the equations $\chi(S^{n,m}_{g,b})
  = \chi(S^{n+2,m}_{g-1,b}) = \chi(S^{n,m+2}_{g,b-1}) =
  \chi(S^{n,m+2}_{g-1,b+1}) = \chi(S^{n+1,m}_{g,b-1}) =
  \chi(S^{n_1,m_1}_{g_1,b_1}) + \chi(S^{n_2,m_2}_{g_2,b_2})$, while
  $p(S^{n,m}_{g,b}) = p(S^{n+2,m}_{g-1,b}) + 1 = p(S^{n,m+2}_{g,b-1})
  + 1 \linebreak[1] = \linebreak[0] p(S^{n,m+2}_{g-1,b+1})
  \linebreak[1] + 1 \linebreak[1] = \linebreak[0] p(S^{n+1,m}_{g,b-1})
  + 1/2 = p(S^{n_1,m_1}_{g_1,b_1}) + p(S^{n_2,m_2}_{g_2,b_2})$ to
  obtain \eqref{MQME1}. Here in both sequences of equations, the last
  one is considered under the assumption of the summation in
  \eqref{premaster}.
\end{sloppypar}
\end{proof}

\begin{rem}
  According to \cite{MMS}, solutions to the quantum master equation in
  geometry are in bijection with the set of wheeled representations of
  a certain wheeled PROP. It would be interesting to see which modular
  operad or, more generally, wheeled PROP would be responsible for the
  above dg BV algebra structure.
\end{rem}

\section{An Algebraic Counterpart}
\label{algebraic}

First of all, consider a linear version of the BV structure of
Section~\ref{BV}. Suppose we have a pair of complexes $H_c$ and $H_o$
of $\C$-vector spaces, whose physical meaning is the state spaces,
including ghosts, of the closed and the open string, respectively.
Suppose these spaces are provided with symmetric bilinear forms
\begin{eqnarray*}
H_c \otimes H_c &\to& \C[1],\\
a\otimes b & \mapsto & (a,b) \in \C,
\end{eqnarray*}
for $a, b$ in $H_c$,
\begin{eqnarray*}
H_o[1] \otimes H_o[1] &\to& \C,\\
a\otimes b & \mapsto & (a,b)' \in \C,
\end{eqnarray*}
for $a, b$ in $H_o$, and
\begin{eqnarray*}
H_o \otimes H_o &\to& \C,\\
a\otimes b & \mapsto & (a,b)'' \in \C,
\end{eqnarray*}
for $a, b$ in $H_o$, which are assumed to be morphisms of complexes.
Suppose also a morphism
\[
\delco: H_c \to \C
\]
of complexes is given.  Form the following space
\[
A := S(H_c) \otimes S(C^\lambda (H_o)),
\]
where $S$ stands for graded symmetric algebra and $C^\lambda$ denotes
the \emph{reduced cyclic complex} (considered without the standard
differential $b$, with $\lambda$ having nothing to do with the string
coupling constant):
\[
C^\lambda (H_o) := \bigoplus_{n=-1}^\infty H_o[1]^{\otimes (n+1)}[-1]
/ (1-t),
\]
where the grading shifts result in placing $a_0 \otimes a_1 \otimes
\dots \otimes a_n$ in degree $\abs{a_0} + \abs{a_1} + \dots +
\abs{a_n} - n$, the operator $t$ is the cyclic-permutation generator:
\[
t(a_0 \otimes a_1 \otimes \dots \otimes a_n) :=
(-1)^{(\abs{a_n}-1)(\abs{a_0} + \dots + \abs{a_{n-1}} - n)} a_n
\otimes a_0 \otimes \dots \otimes a_{n-1}
\]
for any $a_0, a_1, \dots, a_n \in H_o$, and the quotient is taken with
respect to the image of the operator $1-t$, i.e., $C^\lambda$ is the
space of coinvariants of $t$.

The graded vector space $A$ may be endowed with the following
structure of a dg BV algebra. The differential $d$ is just the
internal differential, i.e., the one coming from the differentials on
$H_c$ and $H_o$. The BV operator $\Delta$ is the sum of three components:
\[
\Delta := \Delta_c + \Delta_o + \delco,
\]
which are defined as follows:
\[
\Delta_c (a_1 \otimes a_2 \otimes \dots \otimes a_n) := \sum_{1 \le
  i<j \le n} (-1)^\epsilon a_1 \otimes \dots \otimes \hat{a}_i \otimes
\dots \otimes \hat{a}_j \otimes \dots \otimes a_n (a_i, a_j)
\]
for $a_1, a_2, \dots, a_n \in H_c$ with $(-1)^\epsilon$ being the sign
coming from taking $a_i$ and $a_j$ over to the back of the tensor
product; then $\Delta_c$ is extended by $S(C^\lambda (H_o))$-linearity
to the whole algebra $A$.
\[
\Delta_o := \Delta'_o + \Delta''_o
\]
with the following components:
\begin{multline*}
\Delta'_o (a_0 \otimes a_1 \otimes \dots \otimes a_n)
\\ := \sum_{0 \le
  i<j \le n} (-1)^\epsilon [a_{j+1} \otimes \dots \otimes a_n \otimes
a_0 \otimes \dots \otimes a_{i-1}] \otimes (a_i, a_j)' [a_{i+1} \otimes
\dots \otimes a_{j-1}]\\
\in S^2 (C^\lambda (H_o))
\end{multline*}
for a representative $a_0 \otimes a_1 \otimes \dots \otimes a_n$ of a
cyclic chain in $C^\lambda (H_o)$ with $(-1)^\epsilon$ being the sign
coming from permuting $a_0 \otimes \dots \otimes a_n$ cyclically to
$a_j \otimes \dots \otimes a_n \otimes a_0 \otimes \dots \otimes
a_{j-1}$ and taking $a_j \in H_o[1]$ over to the place right after
$a_i$. The operator $\Delta'_o$ is extended to the symmetric algebra
$S(C^\lambda(H_o))$ as a graded derivation and to the whole tensor
product $A$ by $S(H_c)$-linearity.
\begin{multline*}
\Delta''_o (A_1 \otimes A_2 \otimes \dots \otimes A_n)\\
:= \sum_{1 \le
  i<j \le n} (-1)^\epsilon A_1 \otimes \dots \otimes \hat{A}_i \otimes
\dots \otimes \hat{A}_j \otimes \dots \otimes A_n \otimes \{A_i, A_j\}'',
\end{multline*}
where $A_1, \dots, A_n \in C^\lambda (H_o)$ and $(-1)^\epsilon$ is the
sign coming from taking $A_i$ and $A_j$ to the right in $S(C^\lambda
(H_o))$, and
\begin{multline}
\label{o-bracket}
  \{a_0 \otimes \dots \otimes a_k , b_0 \otimes \dots \otimes b_l\}''\\
  := \sum_{\substack{0 \le i \le k\\0 \le j \le l}} (-1)^\epsilon
  a_{i+1} \otimes \dots \otimes a_k \otimes a_0 \otimes \dots \otimes
  a_{i-1} \otimes (a_i, b_j)'' b_{j+1} \otimes \dots \otimes b_l
  \otimes b_0 \otimes \dots \otimes b_{j-1}
\end{multline}
for $[a_0 \otimes \dots \otimes a_k]$ and $[b_0 \otimes \dots \otimes
b_l] \in C^\lambda (H_o)$ and $(-1)^\epsilon$ being the sign coming
from permuting $a_0 \otimes \dots \otimes a_k$ into $a_{i+1} \otimes
\dots \otimes a_k \otimes a_0 \otimes \dots \otimes a_{i-1} \otimes
a_i$ and $b_0 \otimes \dots \otimes b_l$ into $b_j \otimes b_{j+1}
\otimes \dots \otimes b_l \otimes b_0 \otimes \dots \otimes b_{j-1}$
in $C^\lambda (H_o)$. The operator $\Delta''_o$ is extended to $A$ by
$S(H_c)$-linearity.

Finally, let
\[
\delco: S(H_c) \otimes S(C^\lambda (H_o)) \to S(H_c) \otimes
S(C^\lambda (H_o))[1]
\]
be the extension of the map $\delco: H_c \to \C \subset C^\lambda
(H_o)[1] \subset S(C^\lambda (H_o))[1]$, as before, as an $S(C^\lambda
(H_o))$-linear derivation.

\begin{prop}
  Together with the standard graded commutative multiplication, the
  operator $\Delta$ defines the structure of a dg BV algebra on the
  complex $A = S(H_c) \otimes S(C^\lambda (H_o))$. If we consider the
  formal power series $A[[\lambda, \sqrt{\hslash}]]$ and modify
  $\Delta$ to be $\Delta:= \Delta_c + \Delta_o + \sqrt{\hslash}
  \delco$, we will get the structure of a dg BV algebra on
  $A[[\lambda, \sqrt{\hslash}]]$ over $\C[[\lambda, \sqrt{\hslash}]]$.
\end{prop}

\begin{rem}
  Note that the ``closed-sector'' component $\Delta_c$ is just another
  form of the standard odd Laplacian $\sum_i \del^2 / \del x^i \del
  \xi^i$, see \cite{schwarz:BV}, in the case when the form $(a,b)$
  defines the structure of an odd symplectic manifold on $H_c^*$ with
  Darboux coordinates $(x^i,\xi^i)$. Also, the antibracket
  \eqref{o-bracket} on the reduced cyclic complex is the one coming
  from the Gerstenhaber bracket on the Hochschild cochain complex,
  while the antibracket induced by $\Delta'_o$ is zero.
\end{rem}

\section{Symmetric Open-Closed TCFTs}
\label{TCFT}

\begin{df}
\label{soctcft}
A \emph{symmetric OC TCFT} is a morphism of dg BV algebras $\phi:
V[[\lambda, \sqrt{\hslash}]] \to A[[\lambda, \sqrt{\hslash}]]$ over
$\C[[\lambda, \sqrt{\hslash}]]$, where $V[[\lambda, \sqrt{\hslash}]]$
is the one from Sections \ref{BV} and \ref{section:QME}, based on
chains in the OC moduli spaces, and $A[[\lambda, \sqrt{\hslash}]]$ is
the algebraic one from Section~\ref{algebraic}, based on the state
spaces $H_c$ and $H_o$.  The morphism must also respect some extra
gradings: the component with $n$ interior punctures and $m_1, \dots,
m_b$ punctures located on $b$ boundary components must map to the
component of $S(H_c) \otimes S(C^\lambda (H_o)) [[\lambda,
\sqrt{\hslash}]]$ with the number of factors from $H_c$ equal to $n$,
the number of factors from $C^\lambda(H_o)$ equal to $b$, and the
numbers of those from $H_o$ in these factors equal to $m_1, \dots,
m_b$, irrespective of the order.
\end{df}

\begin{rem}
  Here the word ``symmetric'' refers to the fact that we symmetrize
  over all inputs and outputs (not mixing up open and closed strings,
  though). Mathematically, this means that we consider the moduli
  space with unlabeled interior and boundary punctures, as in
  Section~\ref{OCM}. It is quite obvious how to ``unsymmetrize'' the
  above definition, which will then turn into a variation of the
  standard one, see \cite{costello:tcft,zwiebach:98}. The idea of
  Definition~\ref{soctcft} is that a map $\phi$ respecting the dot
  products behaves well under disjoint union of Riemann surfaces,
  while the condition of respecting the BV operators means that the
  map $\phi$ behaves well under attaching Riemann surfaces at
  punctures, i.e., the correlators of the theory satisfy a
  \emph{factorization axiom}, as well as $\phi$ behaves well with
  respect to turning interior punctures into degenerated boundary
  components, some sort of open-closed string interaction. The idea
  that an OC TCFT gives rise to a morphism from a BV algebra similar
  to $V$ to another BV algebra is due to Zwiebach, see
  \cite{zwiebach:98}.
\end{rem}

Thus, by definition, the dg BV algebra $V$ constructed in
Section~\ref{BV} is the \emph{universal symmetric OC TCFT}. In
particular, if $S \in V[[\lambda, \sqrt{\hslash}]]$ is a solution to
the QME, then its image $\phi(S) \in A[[\lambda, \sqrt{\hslash}]]$ is
a solution to the QME as well.  On the other hand, the dg BV algebra
$V$ is a particular case of the dg BV algebra of the sigma model, when
the target is a one-point space.

It is widely expected, see \cite{costello:gp,fooo,sullivan:sigma},
that the OC sigma model (A model or Gromov-Witten theory) produces a
Gromov-Witten potential satisfying QME. We summarize this expectation
in the language of symmetric OC TCFT as follows.
\begin{conj}
Sigma model produces an example of a symmetric OC TCFT.
\end{conj}

\section{Algebraic Structures}
\label{sh}

Suppose we have a symmetric OC TCFT $V[[\lambda, \sqrt{\hslash}]] \to
A[[\lambda, \sqrt{\hslash}]]$, as in the previous section. We claim
that this data implies certain algebraic structures on $H_c$ and
$H_o$, as well as on the space $U$ generated by the geometric cycles
in the connected OC moduli spaces.

\subsection{The $L_\infty$ structure coming from spheres}

Let us consider the Riemann sphere with punctures and the
corresponding part of the solution \eqref{solution} to the QME:
\[
S_c := S^{\bullet,0}_{0,0} : = \sum_{n \ge 3} S_{0,0}^{n,0}
\lambda^{2n-4} \hslash^{n/2-1}.
\]
Since the image of operator $\Delta_c$ corresponds to strictly
positive genus and the image of $\delco$ to strictly positive number
of boundary components, the QME rewrites in this case as
\begin{equation}
\label{cCME}
dS_c + \frac{1}{2} \{ S_c, S_c\} = 0,
\end{equation}
which is a \emph{Maurer-Cartan equation}, also known as the
\emph{Classical Master Equation $($CME$)$}. This equation implies that
$d + \{S_c, -\}$ is a differential of degree one on the graded
commutative algebra $V[[\lambda, \sqrt{\hslash}]]$.

If we have a symmetric OC TCFT $\phi$, then, since $\phi$ preserves
grading by $b$, the image $\phi(S_c)$ of $S_c$, whose $b$-grading is
zero, must lie in $S(H_c)[[\lambda, \sqrt{\hslash}]]$ and the operator
$d + \{\phi(S_c),-\}$ is a differential on $S (H_c) [[\lambda,
\sqrt{\hslash}]]$. This translates into the structure of an $L_\infty$
coalgebra on $H_c [1]$, according to the following definition.

\begin{df}
  The structure of an $L_\infty$ \emph{coalgebra} on a dg vector space
  $\g$ is a (continuous) degree-one differential $D$ on the completed
  graded symmetric algebra $\hat S(\g[-1])$, the completion being
  taken with respect to the $(\g[-1])$-adic topology.
\end{df}

\begin{rem}
  Since such a derivation is determined by its value on the subspace
  $\g[-1]$ of generators of the symmetric algebra, this structure
  gives rise to a collection of linear maps $D_k: \g[-1] \to
  \g^{\otimes k} [-k]$, $k \ge 1$, interpreted as higher cobrackets,
  satisfying identities dual to those satisfied by brackets in an
  $L_\infty$ algebra. These identities come from breaking the equation
  $D^2 = 0$ into components $D_k$. In the finite-dimensional case,
  $\dim \g < \infty$, this structure is equivalent to the structure of
  an $L_\infty$ algebra on the dual space $\g^*$.
\end{rem}

\begin{prop}
  In a symmetric OC TCFT, the space $H_c[1]$ carries the structure of
  an $L_\infty$ coalgebra.
\end{prop}

\begin{proof}
  The degree-one differential $D(\lambda, \hslash) = d +
  \{\phi(S_c),-\}$ on $S (H_c) [[\lambda, \sqrt{\hslash}]]$ expands as
  follows:
\[
D(\lambda, \hslash) = d + \sum_{n \ge 3} \{\phi( S_{0,0}^{n,0}),-\}
\lambda^{2n-4} \hslash^{n/2-1} = \sum_{k \ge 1} D_k \lambda^{2(k-1)}
\hslash^{(k-1)/2},
\]
where $D_1 := d$ and $D_k := \{ \phi ( S_{0,0}^{k+1,0}),-\}$ for $k
\ge 2$ are operators $D_k : H_c \to S^k (H_c)$. The fact that
$D(\lambda, \hslash)$ is a differential on $S(H_c)[[\lambda,
\sqrt{\hslash}]]$ implies that $D = D(1,1)$ is a differential on the
completed symmetric algebra $\hat S(H_c)$ and thereby defines the
structure of an $L_\infty$ coalgebra on $H_c [1]$.
\end{proof}

Note also that the relative homology classes of the fundamental chains
in $\MUmnv{k+1}{0}{0,0}$ generate the $L_\infty$ operad, while their
boundaries are generated by twist-attachments of lower-dimensional
fundamental chains, producing the defining relations of the $L_\infty$
operad, see \cite{ksv}.  These classes are stable with respect to
permutations of punctures, and their symmetrizations satisfy the CME
\eqref{cCME}, therefore, the $L_\infty$ coalgebra structure on
$H_c[1]$ coincides, up to duality, with the $L_\infty$ algebra
structure on $H_c$ constructed in \cite{zwiebach:93}, see \cite{ksv}.

\subsection{The cyclic $A_\infty$ structure coming from disks}

Kontsevich introduced in \cite{kontsevich:feynman} the notion of a
cyclic $A_\infty$ algebra as an $A_\infty$ algebra with an invariant
inner product. This notion has different variations, called symplectic
$A_\infty$ algebras \cite{kontsevich:formal, hamilton-lazarev:char},
Calabi-Yau $A_\infty$ algebras \cite{kontsevich:Utalk,costello:tcft},
and Sullivan-Wilson's homotopy open Frobenius algebras. The following
definition is motivated by deformation theory of cyclic $A_\infty$
algebras of \cite{penkava-schwarz}.
\begin{df}
  A \emph{cyclic $A_\infty$ algebra} structure on a $($dg$)$ vector
  space $H$ with a symmetric bilinear form $H \otimes H \to \C$ is a
  solution $M = \sum_{k \ge 2} m_{k+1} \lambda^{k-1}$ for $m_{k+1} \in
  H[1]^{\otimes (k+1)}[-1]/(1-t) \subset C^\lambda (H)$ to the CME
  \[
  dM + \frac{1}{2}\{M,M\} = 0
  \]
  in the reduced cyclic chain complex $C^\lambda (H)[[\lambda]]$
  provided with the antibracket given by \eqref{o-bracket} and the
  internal differential $d$ coming from that on $H$.
\end{df}
Before describing this structure in the presence of an OC TCFT, let us
show that this definition is equivalent to the notion of an $A_\infty$
algebra with an invariant inner product, under the assumption that the
underlying graded vector space has finite-dimensional graded
components and the inner product $(,): H \otimes H \to \C$ is
nondegenerate componentwise. A straightforward computation shows that
the CME $dM + \frac{1}{2}\{M,M\} = 0$ is equivalent to the equation
\begin{equation}
\label{hat-CME}
d \hat M + \frac{1}{2}\{\hat M,\hat M\} = 0,
\end{equation}
where $\hat M = \sum \hat m_k$ with $\hat m_k \in \Hom (H^{\otimes k},
H)$ being the operator defined by
\begin{equation}
\label{adjoint}
\hat m_k (v_1 \otimes \dots \otimes v_k) := \sum_{i=0}^k (-1)^\epsilon
m^{(i)}_{k+1} (m^{(i+1)}_{k+1}, v_1) \dots (m^{(i+k)}_{k+1}, v_k),
\end{equation}
where $m_{k+1} = m_{k+1}^{(0)} \otimes m_{k+1}^{(1)} \otimes \dots
\otimes m_{k+1}^{(k)} \in H^{\otimes (k+1)}[k]$ in the standard
notation skipping the summation sign and $(-1)^\epsilon$ is the sign
coming from permuting $m_{k+1}^{(0)} \otimes m_{k+1}^{(1)} \otimes
\dots \otimes m_{k+1}^{(k)} \otimes v_1 \otimes \dots \otimes v_k$
into $m^{(i)}_{k+1} \otimes m^{(i+1)}_{k+1} \otimes v_1 \otimes \dots
\otimes m^{(i+k)}_{k+1} \otimes v_k$. The equation $ d \hat M +
\frac{1}{2}\{\hat M,\hat M\} = 0$ takes place in the Hochschild
complex $\Hom (H[1]^{\otimes \bullet}, H[1])$, with just the internal
differential $d$ and the antibracket being the Gerstenhaber bracket.
If we think of that differential as an element $d \in \Hom (H[1],
H[1])$, then the CME \eqref{hat-CME} is equivalent to $\{ d + \hat M,
d + \hat M\} = 0$, known to be equivalent to the fact that the formal
series $d + \hat M$ defines an $A_\infty$ structure on $H$, as the
Gerstenhaber bracket on the Hochschild complex is the same as the
commutator of Hochschild cochains identified with derivations of the
tensor coalgebra on $H$, \cite{stasheff:g-bracket}.  Equation
\eqref{adjoint} makes it obvious that the higher products $\hat m_k$
are invariant with respect to the inner product.

Now, if we take the part of our solution $S$ to the QME corresponding
to the disk with boundary punctures,
\[
S_o := S^{0,\bullet}_{0,1} : = \sum_{m \ge 3} S_{0,1}^{0,m}
\lambda^{m-2},
\]
it will satisfy the CME
\begin{equation*}
dS_o + \frac{1}{2} \{ S_o, S_o\} = 0.
\end{equation*}
If $\phi: V[[\lambda,\sqrt{\hslash}]] \to A[[\lambda,\sqrt{\hslash}]]$
is an OC TCFT, then by grading arguments, the image $\phi(S_o)$ will
be contained in $C^\lambda (H_o)[[\lambda]]$ and also satisfy the
CME. This implies the following result.

\begin{prop}
  In a symmetric OC TCFT, the space $H_o$ carries the structure of a
  cyclic $A_\infty$ algebra.
\end{prop}

\begin{rem}
  Note that the solution \eqref{solution} of the QME modulo
  $\sqrt{\hslash}$ takes into account only the terms with $p-\chi =
  0$, which implies $g=n=0$, $b=1$. Those terms correspond to the disk
  with punctures on the boundary, and one can think of $S_o$ as $S
  \mod \sqrt{\hslash}$. The apperance of an $A_\infty$ structure is
  geometrically due to the fact that the spaces $\MUmnv{0}{m}{0,1}$
  are the homeomorphic to the Stasheff associahedra $K_{m-1}$,
  \cite{kontsevich:feynman}.
\end{rem}

\bibliographystyle{amsalpha}
\bibliography{string}

\end{document}